\documentclass{gtart_a}
\pdfoutput=1
\usepackage{pictexwd}

\title[Four-dimensional symplectic cobordisms containing three-handles]{Four-dimensional symplectic cobordisms\\containing three-handles}
\author{David T Gay}
\givenname{David T}
\surname{Gay}
\address{Department of Mathematics and Applied 
Mathematics\\University of Cape Town\\\newline
Private Bag X3\\Rondebosch 7701\\South Africa}
\email{dgay@maths.uct.ac.za}
\urladdr{}

\volumenumber{10}
\issuenumber{}
\publicationyear{2006}
\papernumber{42}
\lognumber{0765}
\startpage{1749}
\endpage{1759}

\doi{}
\MR{}
\Zbl{}


\received{22 June 2006}
\revised{}
\accepted{13 October 2006}
\published{28 October 2006}
\publishedonline{28 October 2006}
\proposed{Yasha Eliashberg}
\seconded{Peter Ozsv\'ath, Tomasz Mrowka}
\corresponding{}
\editor{CPR}
\version{}

\arxivreference{math.GT/0606402}

\AtBeginDocument{\let\bar\wbar\let\tilde\wtilde}

\makeatletter
\def\cnewtheorem#1[#2]#3{\newtheorem{#1}{#3}[section]
\expandafter\let\csname c@#1\endcsname\c@theorem}

\newtheorem{theorem}{Theorem}

\cnewtheorem{proposition}[theorem]{Proposition}
\cnewtheorem{corollary}[theorem]{Corollary}
\cnewtheorem{conjecture}[theorem]{Conjecture}
\cnewtheorem{addendum}[theorem]{Addendum}

\theoremstyle{definition}
\cnewtheorem{definition}[theorem]{Definition}
\cnewtheorem{example}[theorem]{Example}
\cnewtheorem{xca}[theorem]{Exercise}
\cnewtheorem{remark}[theorem]{Remark}

\cnewtheorem{notation}[theorem]{Notation}
\cnewtheorem{question}[theorem]{Question}

\makeatother  

\numberwithin{equation}{section}

\font\thinlinefont=cmr5

\def\N{\mathbb N}

\begin{document}

\begin{abstract}
  We construct four-dimensional symplectic cobordisms between contact
  three-man\-ifolds generalizing an example of Eliashberg. One key
  feature is that any handlebody decomposition of one of these
  cobordisms must involve three-handles. The other key feature is that
  these cobordisms contain chains of symplectically embedded
  two-spheres of square zero. This, together with standard gauge
  theory, is used to show that any contact three-manifold of non-zero
  torsion (in the sense of Giroux) cannot be strongly symplectically
  fillable. John Etnyre pointed out to the author that the same
  argument together with compactness results for pseudo-holomorphic
  curves implies that any contact three-manifold of non-zero torsion
  satisfies the Weinstein conjecture.  We also get examples of weakly
  symplectically fillable contact three-manifolds which are (strongly)
  symplectically cobordant to overtwisted contact three-manifolds,
  shedding new light on the structure of the set of contact
  three-manifolds equipped with the strong symplectic cobordism
  partial order.
\end{abstract}

\maketitle

\section{Main Theorem}

Eliashberg~\cite{EliashT3} showed that certain contact structures
$\xi_n$ on $T^3$, although weakly symplectically fillable, are not
strongly symplectically fillable (for $n>1$), by constructing a
symplectic cobordism from $(T^3,\xi_n)$ to the disjoint union of $n$
copies of the standard contact $S^3$. An interesting feature, then, is
that this is a cobordism from a connected contact $3$--manifold to a
disconnected contact $3$--manifold and thus necessarily contains a
$3$--handle. By contrast, most other constructions of symplectic
cobordisms from contact $3$--manifolds to contact $3$--manifolds (as
in Etnyre and Honda~\cite{EtnyreHonda}) are built out of elementary
$0$--, $1$-- and $2$--handle cobordisms as in
Eliashberg~\cite{EliashHandles} or Weinstein~\cite{WeinsteinHandles}. To
this author's knowledge there is no model for an elementary
contact-to-contact $3$--handle symplectic cobordism. Thus it would be
interesting to isolate the $3$--handle from Eliashberg's
$T^3$--to--$S^3$ cobordisms and then generalize Eliashberg's
nonfillability results. We have not succeeded in usefully isolating
the $3$--handle but we have localized the construction to a certain
extent, with new and interesting consequences for nonfillability and
for the cobordism relation in general. (For background on varieties of
fillability and known results, see Etnyre and
Honda~\cite{EtnyreConcave,EtnyreHonda}.)

Given any real numbers $a > b$, let $M_{a,b} = [a,b] \times S^1 \times
S^1$, with coordinates $(t,x,y)$, and consider the contact structure
$\xi_{a,b} = \ker (\cos (\pi t) dy + \sin(\pi t) dx)$. Given any $c >
0$ let $(M_c,\xi_c) = (M_{0,c},\xi_{0,c})$.  Also, using ``toric
coordinates'' $(p,q)$ on $\R^2$ with $p = r^2/2$ and $q = \theta$, and
given any $a>0$, let $T_a$ be the solid torus $S^1 \times \{p \leq
a\}$ (with $\alpha$ the $S^1$--coordinate), and let $\eta_a = \ker (\cos(\pi
p) d\alpha + \sin(\pi p) dq)$.

Our main result is:
\begin{theorem} \label{TCobordism} For any $k \in \N$, with $k > 1$,
  there exists a smooth cobordism-with-sides $X_k$ from $M_{k+1/2}$ to
  $T_{1/2} \amalg T_{1/2}$ containing a chain of embedded $2$--spheres
  $S_1, \ldots, S_{2k-2}$, with a symplectic form $\omega_k$,
  satisfying the following properties:
\begin{itemize}
\item Each $S_i$ is symplectic and all intersections between $S_i$'s
  are positive.
\item $S_i \cdot S_j = 0$ unless $j-i = \pm 1$, in which case $S_i
  \cdot S_{i+1} = +1$.
\item There is a nowhere-zero Liouville vector field $V$ defined in a
  neighborhood of the bottom and sides, pointing in along the bottom
  and parallel to the sides, inducing the contact structure
  $\xi_{k+1/2}$ on the bottom $M_{k+1/2}$ (i.e. $\xi_{k+1/2} = \ker
  (\imath_V \omega |_{M_{k+1/2}})$).
\item Along the top, the contact structure $\xi_{1/2}$ on each copy of
  $T_{1/2}$ is dominated by $\omega_k$ (i.e. $\omega_k |_{\xi_{1/2}} >
  0$), and agrees with $\ker (\imath_V \omega|_{T_{1/2}})$ on a
  neighborhood of $\partial T_{1/2}$.
\end{itemize}
\end{theorem}

A nearly identical theorem could be stated producing a chain of
(fewer) spheres of self-intersection $+1$. The version above is
slightly more convenient to state and prove because a square (the
moment map image of $S^2 \times S^2$) is slightly easier to work with
in cartesian coordinates than a triangle (the moment map image of $\C
P^2$). We have not stated which curves on $\partial M_{k+1/2}$
correspond to meridians in $T_{1/2} \amalg T_{1/2}$ as we do not need
this information for our applications; if needed, this information is
extractable from the proof.  Also, for the applications in this paper
we only need the case $k=2$; presenting only this case would shorten
the proof, but we believe that the general result is interesting in
its own right so we present it in full.

\section{Corollaries}

Before the proof, we discuss some corollaries. As a prerequisite we
have:
\begin{definition}[Giroux~\cite{GirouxMeme,GirouxInfinite,GirouxBifurcations}] Given an isotopy class $C$ of tori in a contact
$3$--manifold $(M,\xi)$, the {\em torsion} of $C$ is the largest
integer $n$ such that $(M_{2n},\xi_{2n})$ can be contactomorphically
embedded in $(M,\xi)$ as a neighborhood of a torus in $C$.
\end{definition}
We will be a bit sloppy, and speak of the torsion of a torus, meaning
the torsion of its isotopy class. 

\begin{corollary} \label{CTorsion} If a contact $3$--manifold
  $(M,\xi)$  contains a torus of torsion greater than or equal
  to $1$ then $(M,\xi)$ is not symplectically fillable.
\end{corollary}

This was apparently a conjecture of Eliashberg's.
It should be noted that Ding and Geiges~\cite{DingGeiges} generalized
Eliashberg's result on nonfillability of contact structures on $T^3$
to contact structures on general $T^2$ bundles over $S^1$, and thus
proved this corollary in certain very special
cases. Ghiggini~\cite{Ghiggini} then proved this for certain Seifert
fibred $3$--manifolds, and most recently Lisca and
Stipsicz~\cite{LiscaStipsicz} proved it for a larger class of
$3$--manifolds, characterized in terms of their Ozsvath--Szabo
invariants.  

\begin{proof}
  Note that a neighborhood of $(M_k,\xi_k)$ always contains
  $(M_{k+\epsilon},\xi_{k+\epsilon})$ for small $\epsilon > 0$, and
  thus that a neighborhood of $(M_2,\xi_2)$ contains a copy of
  $(M_{2+1/q},\xi_{2+1/q})$ for large enough integer $q$. However, for
  any $q \in \Z$, $(M_{2+1/q},\xi_{2+1/q})$ is contactomorphic to
  $(M_{5/2},\xi_{5/2})$ via a contactomorphism of the form $(t,(x,y))
  \mapsto (f(t),L(x,y))$ where $f$ is smooth and increasing and $L$ is
  linear. Thus we now assume that $(M,\xi)$ contains
  $(M_{5/2},\xi_{5/2})$.
  
  Suppose $(M,\xi)$ is the strongly convex boundary of a compact
  symplectic $4$--manifold $(X,\omega)$. Then we can attach the
  cobordism $(X_2,\omega_2)$ along $(M_{5/2},\xi_{5/2})$, attaching a
  trivial cobordism made from the symplectization of $\xi$ along the
  rest of $M$, to produce a symplectic $4$--manifold $(X',\omega')$
  with (weakly) convex boundary containing a symplectic ``hyperbolic
  pair'' $(S_{1}, S_{2})$, with $S_1 \cdot S_2 = +1 = |S_1 \cap S_2|$
  and $S_1 \cdot S_1 = S_2 \cdot S_2 = 0$.
  
  Now cap off $(X',\omega')$ with a concave filling (see
  Etnyre~\cite{EtnyreConcave} or Eliashberg~\cite{EliashConcave}) to
  get a closed symplectic $4$--manifold. Such a concave filling can
  always be constructed so that $b_2^+$ of the filling is
  positive. This is because, if we use the construction
  in~\cite{EtnyreConcave}, the first step is to construct a cobordism
  up to a homology $3$--sphere (Lemma~3.1 of~\cite{EtnyreConcave}), at
  which point the weakly convex boundary can be made into a strongly
  convex boundary (Ohta and Ono~\cite{OhtaOno}), which can then be
  capped off with a concave filling as in Gay~\cite{GayConcave}, in
  which one can explicitly see a surface of positive
  self-intersection. Thus we have a closed symplectic $4$--manifold
  with $b_2^+ > 1$, containing a hyperbolic pair of spheres, which is
  well-known to be impossible: Since the boundary of a neighborhood of
  a hyperbolic pair (with geometric intersection equal to algebraic
  intersection) is $S^3$, this means that the $4$--manifold splits as
  a connected sum of two $4$--manifolds each with $b_2^+ \geq 1$,
  which cannot happen because symplectic $4$--manifolds have nontrival
  Seiberg--Witten invariants (Tuabes~\cite{Taubes}) while connected
  sums of $4$--manifolds with $b_2^+ \geq 1$ have trivial
  Seiberg--Witten invariants (see~\cite{OhtaOno}, for example).
\end{proof}

The next corollary and the proof presented here were explained to the
author by John Etnyre. Recall that the generalized Weinstein
conjecture for a contact $3$--manifold $(M,\xi)$ states that any Reeb
vector field for $\xi$ has a closed orbit.

\begin{corollary}[Etnyre]
Any contact $3$--manifold $(M,\xi)$ containing a torus of non-zero torsion
satisfies the Weinstein conjecture.
\end{corollary}

Note that this is already known to be true on a large class of
non-zero torsion contact $3$--manifolds as a result of the computations
in Bourgeois and Colin~\cite{BourgeoisColin}. 

\begin{proof}[Sketch of proof]
As in the preceding proof, we can construct a concave cap for
$(M,\xi)$, namely a symplectic $4$--manifold $(X,\omega)$ with
strongly concave boundary $(M,\xi)$, containing a symplectic sphere
$S$ of square zero. Attaching the negative symplectization of a
particular contact form $\alpha$ for $\xi$ gives a noncompact manifold
with negative cylindrical end as in Bourgeois, Eliashberg, Hofer,
Wysocki and Zehnder~\cite{BEHWZ}.  The compactness results
in~\cite{BEHWZ} can be used to extend to this setting the techniques
which McDuff~\cite{McDuff} used to understand closed symplectic
$4$--manifolds containing square zero symplectic $2$--spheres. This
shows that, if $R_\alpha$ has no closed orbits, then the moduli space
of $J$--holomorphic spheres homologous to $S$ (for a suitable $J$) is
a compact manifold, leading to a contradiction.
\end{proof}

We can also reprove an old classic, although as the referee pointed
out this does not really constitute a new proof:

\begin{corollary}[Eliashberg~\cite{EliashTight}, Gromov~\cite{Gromov}]
  \label{CFillableImpliesTight} Any weakly symplectically
  semifillable contact $3$--manifold $(M,\xi)$ is tight.
\end{corollary}

\begin{proof}
  Suppose that $(M,\xi)$ is overtwisted but weakly symplectically
  semifillable.  Let $(X,\omega)$ be a weak symplectic semi-filling of
  $(M,\xi)$. Cap off all the other components of $\partial X$ as
  in~\cite{EtnyreConcave} or~\cite{EliashConcave}. Now apply Lemma~3.1
  of either~\cite{EtnyreConcave} or Stipsicz~\cite{Stipsicz} to attach
  a cobordism to $(X,\omega)$ to arrive at a homology $3$--sphere
  contact boundary $(M',\xi')$ and then use~\cite{OhtaOno} to make
  sure the new $4$--manifold is in fact a strong filling of
  $(M',\xi')$. The Legendrian surgeries involved in Lemma~3.1
  of~\cite{EtnyreConcave} can be arranged so as to avoid an
  overtwisted disk in $(M,\xi)$, and the perturbations of the
  symplectic structure in~\cite{OhtaOno} do not change the contact
  structure, so we can assume that $(M',\xi')$ is again overtwisted,
  but now with a strong symplectic filling.
  
  An overtwisted contact $3$--manifold has isotopy classes of tori of
  torsion greater than $0$ (in fact of arbitrarily large torsion).
  This is because, if we perform arbitrarily many full Lutz twists
  along a transverse knot in $(M',\xi')$, we produce a contact
  structure which is isotopic to $\xi'$, using the fact (see
  Geiges~\cite{GeigesSurvey}) that full Lutz twists do not change
  homotopy classes of contact structures and the fact
  (Eliashberg~\cite{EliashOT}) that homotopic, overtwisted contact
  structures are isotopic.  Therefore $(M',\xi')$ contains a copy of
  $(T_{5/2},\eta_{5/2})$, which contains a copy of
  $(M_2,\xi_2)$. \fullref{CTorsion} then yields our contradiction.
\end{proof}

To motivate our last corollary we recall the main results
of~\cite{EtnyreHonda}. There is a natural partial order $\prec$ on the
set $\mathcal{C}$ of closed (possibly disconnected, possibly empty)
positive contact $3$--manifolds: $(M_0,\xi_0) \prec (M_1,\xi_1)$ if
there exists a compact (possibly disconnected) cobordism $X$ from
$M_0$ to $M_1$ equipped with a symplectic form $\omega$ and a
Liouville vector field $V$ defined on a neighborhood of $M_0 \cup
M_1$, pointing in along $M_0$ and out along $M_1$, such that $\xi_i =
\ker (\imath_V \omega |_{M_i})$. (This is what we mean by a ``strong
symplectic cobordism'' and a basic fact is that this relation is
reflexive and transitive but not symmetric.) Etnyre and
Honda~\cite{EtnyreHonda} showed the following two facts:
\begin{itemize}
\item $(M,\xi) \prec \emptyset$ for every contact $3$--manifold $(M,\xi)$.
\item For every {\em connected} overtwisted contact $3$--manifold $(M_o,\xi_o)$
  and for any other {\em connected} contact $3$--manifold $(M,\xi)$, we have
  $(M_o,\xi_o) \prec (M,\xi)$. 
\end{itemize}
The connectedness assumption in the second point is related to the
absence of models for symplectic $3$--handles. Of course $3$--handles
are involved in any cobordism to $\emptyset$, as in the first point,
but those $3$--handles are really upside down $1$--handles, in the sense
that, as elementary cobordisms, if we make the concave end the
bottom and the convex end the top, then they are $1$--handles.

Now it is reasonable to introduce an equivalence relation on
$\mathcal{C}$, whereby $A \sim B$ if $A \prec B$ and $B \prec A$. Then
the partial order $\prec$ descends to $\mathcal{C}/\!\!\sim$, and one
can begin the study of $(\mathcal{C},\prec)$ with the study of
$(\mathcal{C}/\!\!\sim, \prec)$. The above results together with the
fact that overtwisted contact structures are not fillable mean that
there are at least two distinct elements in $\mathcal{C}/\!\!\sim$,
namely $[\emptyset]$ (precisely the strongly symplectically fillable
contact $3$--manifolds), and the equivalence class containing all
connected overtwisted contact $3$--manifolds, which we will call
$\mathcal{O}_1$.  In addition we have that, for every connected $A \in
\mathcal{C}$, $\mathcal{O}_1 \prec [A] \prec [\emptyset]$.  It is also
immediate that, if $\mathcal{O}_i$ is the equivalence class containing
all contact $3$--manifolds with $i$ components each of which is
overtwisted, then $\mathcal{O}_i \prec \mathcal{O}_{i-1}$, and that,
for any contact $3$--manifold $A$ with $i$ components, $\mathcal{O}_i
\prec [A] \prec [\emptyset]$. Thus we have the following natural
questions:
\begin{itemize}
\item Is $\mathcal{O}_i = \mathcal{O}_{i-1}$?
\item Are there any tight contact $3$--manifolds in $\mathcal{O}_1$?
\end{itemize}

We answer both questions affirmatively:
\begin{corollary} \label{CPartialOrder} For every contact
  $3$--manifold $(M,\xi)$ containing a torus $T$ of
  torsion greater than $1$, there is a strong symplectic cobordism
  from $(M,\xi)$ to an overtwisted contact $3$--manifold $(M',\xi')$.
  If $T$ is separating and $M$ is connected then $M'$ will have two
  components, on each of which $\xi'$ is overtwisted.
\end{corollary}
This answers the first question affirmatively because, as noted in the
proof of \fullref{CFillableImpliesTight}, overtwisted contact
$3$--manifolds contain tori of arbitrarily large torsion, which are
separating because they are boundaries of solid tori. This answers the
second question affirmatively because Giroux~\cite{GirouxMeme,GirouxInfinite,GirouxBifurcations} and Colin~\cite{ColinInfinite} have
constructed many examples of weakly fillable, and hence tight, contact
$3$--manifolds with arbitrarily large torsion.  Thus we know the
minimal and maximal elements, $\mathcal{O} = \mathcal{O}_1$ and
$[\emptyset]$, respectively, in $(\mathcal{C}/\!\!\sim,\prec)$, and we are
left with the obvious question:
\begin{question}
  Are there any elements in $\mathcal{C}/\!\!\sim$ other than $\mathcal{O}$ and
  $[\emptyset]$? 
\end{question}
As far as we know, this remains open, but the obvious candidates are
equivalence classes of tight but not weakly fillable contact
$3$--manifolds and of contact $3$--manifolds containing tori of torsion
exactly $1$.

\begin{proof}[Proof of \fullref{CPartialOrder}]
  By the hypotheses, $(M,\xi)$ contains a copy of
  $(M_{9/2},\xi_{9/2})$. We can find real numbers $a,b$ with $0 < a <
  b < 9/2$ so that $(M_{0,a},\xi_{0,a})$ contains a copy of
  $(M_1,\xi_1)$ in its interior, $(M_{a,b},\xi_{a,b})$ contains a copy
  of $(M_2,\xi_2)$ in its interior, and $(M_{b,9/2},\xi_{b,9/2})$
  contains a copy of $(M_1,\xi_1)$ in its interior. Attach the
  cobordism $(X_2,\omega_2)$ from \fullref{TCobordism} to the
  $(M_2,\xi_2)$ inside $(M_{a,b},\xi_{a,b})$ (extending by the
  symplectization of $\xi$ on the rest of $(M,\xi)$). The contact
  $3$--manifold $(M',\xi')$ on the top of the cobordism is then
  obtained from $(M,\xi)$ by removing the $(M_2,\xi_2)$ and replacing
  with two solid tori. Thus the two remaining copies of $(M_1,\xi_1)$
  each end up bounding solid tori, so that the solid tori each contain
  overtwisted meridinal disks. It is clear that if $C$ is separating,
  then $M'$ will be disconnected.
\end{proof}

\section[Proof of \ref{TCobordism}]{Proof of \fullref{TCobordism}}

  This proof uses the technique of seeing $4$--dimensional symplectic
  topology through $2$--dimensional pictures in moment map images and
  their generalizations, especially as developed by
  Symington~\cite{Symington.4from2}, where a leisurely introduction
  may be found.

  Eliashberg's cobordism in~\cite{EliashT3} is an $n$--fold cyclic
  cover of the complement of a neighborhood of a Lagrangian torus in
  $B^4$ equipped with the standard symplectic form. The standard
  moment map on $\R^4$ has image equal to the first quadrant in
  $\R^2$, and this map restricted to $B^4$ has image equal to a
  right-angled triangle with vertices at $(0,0)$, $(1,0)$ and
  $(0,1)$. Removing a disk from the interior of this triangle
  corresponds to removing a neighborhood of a Lagrangian torus from
  $B^4$, and then the $n$--fold cyclic cover in question corresponds to
  the $n$--fold cyclic cover of this punctured right-angled
  triangle. It is this picture that led the author to the following
  construction:
  
  Consider the standard moment map $\mu\co  S^2 \times S^2 \rightarrow
  \R^2$ for the standard torus action, translated and rescaled so that
  $\mu(S^2 \times S^2) = [-1,1] \times [-1,1]$, and let $\omega$ be
  the corresponding symplectic form on $S^2 \times S^2$.  Recall that
  the preimage of a point in the interior of the square $[-1,1] \times
  [-1,1]$ is a Lagrangian torus, that the preimage of a point in the
  interior of an edge is a circle, that the preimage of a vertex is a
  point, that the preimage of an entire edge is a symplectic sphere of
  self-intersection $0$, and that the oriented intersection of two of
  these spheres (covering two adjacent edges) is $+1$.  Let $\Gamma$
  be the complement in $[-1,1] \times [-1,1]$ of the open disk of
  radius $1/4$ centered at $(0,0)$ and let $X = \mu^{-1}(\Gamma)$. In
  other words, $X$ is the complement of an open neighborhood of a
  Lagrangian torus in $S^2 \times S^2$. We will use coordinates
  $(p_1,p_2)$ on $\R^2$; recall that these are paired with angular
  coordinates $(q_1,q_2)$ on $S^2 \times S^2$ so that $\omega = dp_1
  \wedge dq_1 + dp_2 \wedge dq_2$.

  Now let $\tilde{\Gamma}$ be the infinite cyclic cover of $\Gamma$ (with
  covering map $p\co  \tilde{\Gamma} \rightarrow \Gamma$) and let
  $(\tilde{X},\tilde{\omega})$ be the corresponding infinite cyclic cover of
  $(X,\omega)$ (with covering map $\pi\co  \tilde{X} \rightarrow X$), so that the
  natural moment map $\mu \circ \pi$ for $\tilde{X}$ factors through a
  locally toric fibration $\tilde{\mu}\co  \tilde{X} \rightarrow
  \tilde{\Gamma}$.
  
\begin{figure}
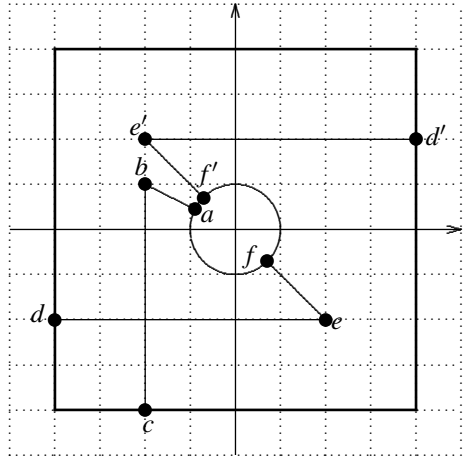

    \begin{center}
\mbox{\beginpicture\small
\setcoordinatesystem units <0.63cm,0.63cm>
\unitlength=0.63cm
\linethickness=1pt
\setplotsymbol ({\makebox(0,0)[l]{\tencirc\symbol{'160}}})
\setshadesymbol ({\thinlinefont .})
\setlinear
%
%
\linethickness=5pt
\setplotsymbol ({\makebox(0,0)[l]{\tencirc\symbol{'164}}})
{\plot  7.620 18.098  7.620 18.098 /
%
%
\linethickness=5pt
\setplotsymbol ({\makebox(0,0)[l]{\tencirc\symbol{'164}}})
\plot 13.335 18.098 13.335 18.098 /
%
%
\linethickness= 0.500pt
\setplotsymbol ({\thinlinefont .})
\plot 12.103 19.327 13.335 18.098 /
}%
%
%
\linethickness=5pt
\setplotsymbol ({\makebox(0,0)[l]{\tencirc\symbol{'164}}})
{\plot 12.099 19.327 12.099 19.327 /
%
%
\linethickness= 0.500pt
\setplotsymbol ({\thinlinefont .})
\putrule from 13.335 18.098 to  7.620 18.098
}%
%
%
\linethickness=5pt
\setplotsymbol ({\makebox(0,0)[l]{\tencirc\symbol{'164}}})
{\plot 15.240 21.907 15.240 21.907 /
%
%
\linethickness=5pt
\setplotsymbol ({\makebox(0,0)[l]{\tencirc\symbol{'164}}})
\plot  9.525 21.907  9.525 21.907 /
%
%
\linethickness= 0.500pt
\setplotsymbol ({\thinlinefont .})
\plot 10.757 20.678  9.525 21.907 /
}%
%
%
\linethickness=5pt
\setplotsymbol ({\makebox(0,0)[l]{\tencirc\symbol{'164}}})
{\plot 10.761 20.678 10.761 20.678 /
%
%
\linethickness= 0.500pt
\setplotsymbol ({\thinlinefont .})
\putrule from  9.525 21.907 to 15.240 21.907
}%
%
%
\linethickness= 0.500pt
\setplotsymbol ({\thinlinefont .})
{\ellipticalarc axes ratio  0.953:0.953  360 degrees 
	from 12.383 20.003 center at 11.430 20.003
}%
%
%
\linethickness=1pt
\setplotsymbol ({\makebox(0,0)[l]{\tencirc\symbol{'160}}})
{\putrectangle corners at  7.620 23.812 and 15.240 16.192
}%
%
%
\linethickness= 0.500pt
\setplotsymbol ({\thinlinefont .})
{\putrule from 11.430 15.240 to 11.430 24.765
%
%
\plot 11.509 24.479 11.430 24.765 11.351 24.479 /
}%
%
%
\linethickness= 0.500pt
\setplotsymbol ({\thinlinefont .})
{\putrule from  6.668 20.003 to 16.192 20.003
%
%
\plot 15.907 19.923 16.192 20.003 15.907 20.082 /
}%
%
%
\linethickness= 0.500pt
\setplotsymbol ({\thinlinefont .})
\setdots < 0.0953cm>
{\plot  7.620 15.240  7.620 24.765 /
}%
%
%
\linethickness= 0.500pt
\setplotsymbol ({\thinlinefont .})
{\plot  8.572 15.240  8.572 24.765 /
}%
%
%
\linethickness= 0.500pt
\setplotsymbol ({\thinlinefont .})
{\plot 10.478 15.240 10.478 24.765 /
}%
%
%
\linethickness= 0.500pt
\setplotsymbol ({\thinlinefont .})
{\plot 12.383 15.240 12.383 24.765 /
}%
%
%
\linethickness= 0.500pt
\setplotsymbol ({\thinlinefont .})
{\plot 13.335 15.240 13.335 24.765 /
}%
%
%
\linethickness= 0.500pt
\setplotsymbol ({\thinlinefont .})
{\plot 14.287 15.240 14.287 24.765 /
}%
%
%
\linethickness= 0.500pt
\setplotsymbol ({\thinlinefont .})
{\plot 15.240 15.240 15.240 24.765 /
}%
%
%
\linethickness= 0.500pt
\setplotsymbol ({\thinlinefont .})
{\plot 16.192 15.240 16.192 24.765 /
}%
%
%
\linethickness= 0.500pt
\setplotsymbol ({\thinlinefont .})
{\plot  6.668 24.765 16.192 24.765 /
}%
%
%
\linethickness= 0.500pt
\setplotsymbol ({\thinlinefont .})
{\plot  6.668 23.812 16.192 23.812 /
}%
%
%
\linethickness= 0.500pt
\setplotsymbol ({\thinlinefont .})
{\plot  6.668 22.860 16.192 22.860 /
}%
%
%
\linethickness= 0.500pt
\setplotsymbol ({\thinlinefont .})
{\plot  6.668 19.050 16.192 19.050 /
}%
%
%
\linethickness= 0.500pt
\setplotsymbol ({\thinlinefont .})
{\plot  6.668 17.145 16.192 17.145 /
}%
%
%
\linethickness= 0.500pt
\setplotsymbol ({\thinlinefont .})
{\plot  6.668 16.192 16.192 16.192 /
}%
%
%
\linethickness= 0.500pt
\setplotsymbol ({\thinlinefont .})
{\plot  6.668 15.240 16.192 15.240 /
}%
%
%
\linethickness=5pt
\setplotsymbol ({\makebox(0,0)[l]{\tencirc\symbol{'164}}})
\setsolid
{\plot  9.525 20.955  9.525 20.955 /
%
%
\linethickness=5pt
\setplotsymbol ({\makebox(0,0)[l]{\tencirc\symbol{'164}}})
\plot  9.525 16.192  9.525 16.192 /
%
%
\linethickness= 0.500pt
\setplotsymbol ({\thinlinefont .})
\putrule from  9.525 20.955 to  9.525 16.192
}%
%
%
\linethickness=5pt
\setplotsymbol ({\makebox(0,0)[l]{\tencirc\symbol{'164}}})
{\plot 10.579 20.428 10.579 20.428 /
%
%
\linethickness= 0.500pt
\setplotsymbol ({\thinlinefont .})
\plot 10.569 20.432  9.525 20.955 /
}%
%
%
\linethickness= 0.500pt
\setplotsymbol ({\thinlinefont .})
\setdots < 0.0953cm>
{\plot  6.668 18.098 16.192 18.098 /
}%
%
%
\linethickness= 0.500pt
\setplotsymbol ({\thinlinefont .})
{\plot  6.668 15.240  6.668 24.765 /
}%
%
%
\linethickness= 0.500pt
\setplotsymbol ({\thinlinefont .})
{\plot  6.668 21.907 16.192 21.907 /
}%
%
%
\linethickness= 0.500pt
\setplotsymbol ({\thinlinefont .})
{\plot  6.668 20.955 16.192 20.955 /
}%
%
%
\linethickness= 0.500pt
\setplotsymbol ({\thinlinefont .})
{\plot  9.525 15.240  9.525 24.765 /
}%
%
%
\put{$b$%
} [lB] <0pt,1pt> at 9.258 21.087
%
%
\put{$c$%
} [lB] <0pt,-2pt> at  9.423 15.83
%
%
\put{$a$%
} [lB] <1pt,0pt> at 10.6 20.16
%
%
\put{$d$%
} [lB] <-2pt,0pt> at  7.176 18.058
%
%
\put{$f$%
} [lB] <-2pt,0pt> at 11.725 19.3
%
%
\put{$e$%
} [lB] at 13.407 17.899
%
%
\put{$d'$%
} [lB] at 15.401 21.797
%
%
\put{$f'$%
} [lB] <0pt,1pt> at 10.647 20.95
%
%
\put{$e'$%
} [lB] at  9.15 22.0
\linethickness=0pt
\putrectangle corners at  6.642 24.790 and 16.218 15.215
\endpicture}
    \caption{The six points in $\Gamma$ used to construct $\chi_k$
      (primes on labels indicate the case when $k$ is odd)}
    \label{FSixPoints}
    \end{center}
  \end{figure}

  Our cobordism $(X_k,\omega_k)$ will be a subset of
  $(\tilde{X},\tilde{\omega})$ which is the preimage via $\tilde{\mu}$ of a
  subset $\chi_k$ of $\tilde{\Gamma}$. We describe $\chi_k$ in the next
  few paragraphs. (See \fullref{FSixPoints}.)

  Consider the following six points in $\Gamma$:
  \begin{itemize}
  \item $a$ is the intersection of the circle of radius $1/4$ centered at
  $(0,0)$ and the line segment from $(0,0)$ to $(-1/2,1/4)$.
  \item $b = (-1/2,1/4)$
  \item $c = (-1/2,-1)$
  \item If $k$ is even then $d = (-1,-1/2)$, otherwise $d = (1,1/2)$.
  \item If $k$ is even then $e = (1/2,-1/2)$, otherwise $e = (-1/2,1/2)$.
  \item If $k$ is even then $f$ is the intersection of the circle of radius
  $1/4$ centered at $(0,0)$ and the line segment from $(0,0)$ to $(1/2,-1/2)$,
  otherwise $f$ is the intersection of the same circle and the line segment
  from $(0,0)$ to $(-1/2,1/2)$.
  \end{itemize}
  
  We will want to work with lifts of these points in $\tilde{\Gamma}$, to
  which end we establish the following conventions: If we use polar
  coordinates $(r,\theta)$ on $\Gamma$, where $\theta \in \R/2\pi \Z$, then we
  can naturally lift to coordinates $(r,\tilde{\theta})$ on $\tilde{\Gamma}$,
  where $\tilde{\theta} \in \R$. Then, for any point $p \in \Gamma$ and any
  integer $i$, $p_i \in \tilde{\Gamma}$ will denote the unique lift of $p$ such
  that $\tilde{\theta}(p_i) \in [2\pi i, 2\pi(i+1))$. Also, when we speak of
  straight line segments in $\tilde{\Gamma}$, we mean arcs that project to
  straight line segments in $\Gamma$.
  
  Let $n = \lceil{k/2}\rceil$. The set $\chi_k \subset \tilde{\Gamma}$ is the
  compact subset of $\tilde{\Gamma}$ bounded by:
  \begin{itemize} 
  \item the straight line segment $A_0$ from $a_0$ to $b_0$,
  \item the straight line segment $B_0$ from $b_0$ to $c_0$,
  \item the part $C$ of the outer (right-angled, polygonal) boundary of
    $\tilde{\Gamma}$ going from $c_0$ to $d_n$,
  \item the straight line segment $B_n$ from $d_n$ to $e_n$,
  \item the straight line segment $A_n$ from $e_n$ to $f_n$, and
  \item the part $D$ of the inner (round) boundary of $\tilde{\Gamma}$ going
    from $f_n$ back to $a_0$.
  \end{itemize}
  
  This defines the set $\chi_k$, and then our cobordism is $X_k =
  \tilde{\mu}^{-1}(\chi_k)$, with $\omega_k =
  \tilde{\omega}|_{X_k}$. The bottom boundary is
  $\tilde{\mu}^{-1}(D)$, the sides are $\tilde{\mu}^{-1}(A_0 \cup
  A_n)$, and the top boundary is $\tilde{\mu}^{-1}(B_0 \cup B_n)$. The
  chain of symplectic spheres is the chain of preimages of the
  straight line segments of length $2$ in $C$. (Note that $C$ starts
  and ends with segments of length $3/2$ and is otherwise composed of
  $2k-2$ segments each of length $2$, each meeting the next at a right
  angle, so that standard toric geometry shows that these spheres
  satisfy all the advertised properties.) Standard toric geometry also
  shows that $\tilde{\mu}^{-1}(D) \cong M_{1-p,k+3/2+q}$, where $2\pi
  p = \tan^{-1}(1/2)$ and $2\pi q = \tan^{-1}(1)$, and that
  $\tilde{\mu}^{-1}(B_0)$ and $\tilde{\mu}^{-1}(B_n)$ are solid tori.
  
  Symington~\cite{Symington.4from2} has observed that, if we pick any point $P
  \in \R^2$ and let $W$ be the outward pointing radial vector field centered
  at $P$, then $P$ lifts to a Liouville vector field $V$ in a toric symplectic
  $4$--manifold $(X,\omega)$ via the moment map $\mu\co  X \rightarrow \R^2$. If
  $E$ is an edge of the polygonal moment map image then $V$ will be defined on
  $\mu^{-1}(E)$ if and only if $V$ is tangent to $E$. (This can be seen 
  by translating the moment map image so that $P = (0,0)$, in which case $V =
  p_1 \partial_{p_1} + p_2 \partial_{p_2}$. These specific coordinates allow
  us also to read off the precise contact structure induced by $V$ on any
  transverse $3$--manifold.)
  
  The Liouville vector field advertised in the theorem is then the
  lift of the radial vector field on $\Gamma$ centered at $(0,0)$,
  i.e. $V = p_1 \partial_{p_1} + p_2 \partial_{p_2}$. This is
  transverse to $\tilde{\mu}^{-1}(D)$ and parallel to
  $\tilde{\mu}^{-1}(B_0)$ and $\tilde{\mu}^{-1}(B_n)$. A direct
  calculation shows that the induced contact structure on
  $\tilde{\mu}^{-1}(D) \cong M_{1-p,k+3/2+q}$ is precisely
  $\xi_{1-p,k+3/2+q}$. Because the angle $\pi(1-k)$ determines the ray
  passing through $(-1,2)$ and the angle $\pi(k+3/2+q)$ determines the
  ray passing through $(1,-1)$, a linear transformation in the $S^1
  \times S^1$ factor shows that this is contactomorphic to
  $(M_{k+1/2}, \xi_{k+1/2})$.

  Now we construct contact structures $\zeta_0$ on
  $\tilde{\mu}^{-1}(B_0)$ and $\zeta_n$ on $\tilde{\mu}^{-1}(B_n)$
  which are dominated by $\omega_k$ and such that both
  $(\tilde{\mu}^{-1}(B_0),\zeta_0)$ and
  $(\tilde{\mu}^{-1}(B_n),\zeta_n)$ are contactomorphic to
  $(T_{1/2},\eta_{1/2})$. To this end, we construct vector fields
  $W_0$ and $W_n$ along $p(B_0)$ and $p(B_n)$ in $\Gamma$ as follows:
  Let $W$ be the radial vector field centered at $(0,0)$, let $U_0$ be
  the radial vector field centered at $(3,-1)$ and let $U_n$ be the
  radial vector field centered at $(-1,3)$ if $n$ is even or at
  $(1,-3)$ if $n$ is odd. Recall that $p(B_0)$ is the line segment
  from $b$ to $c$. We describe $W_0$ as we move from $b$ to $c$.  At
  the beginning, near $b$, $W_0$ is equal to $W$, and then in a short
  interval $W_0$ monotonically interpolates from $W$ to $U_0$, and
  then $W_0$ is equal to $U_0$ on the rest of $p(B_0)$. This can be
  done so that, as we move from $b$ to $c$, $W_0$ is turning
  monotonically counterclockwise relative to the $(p_1,p_2)$
  coordinate system. Similarly, $p(B_n)$ is the line segment from $e$
  to $d$, and $W_n$ interpolates from $W$ to $U_n$ so that $W_n$
  monotonically rotates counterclockwise as we move from $e$ to
  $d$. Now lift $W_0$ and $W_n$ to vector fields $V_0$ and $V_k$ on
  $X_k$ (as we did for the radial vector fields in the preceding
  paragraph, except that now these are not Liouville), and let
  $\zeta_i = \imath_{V_i} \omega_k |_{\tilde{\mu}^{-1}(B_i)}$, for $i
  = 0$ and $n$. Although $V_i$ is not Liouville, nevertheless
  $\zeta_i$ will be contact simply because $W_i$ is transverse to
  $p(B_i)$ and rotates monotonically counterclockwise as we move along
  the length of $p(B_i)$. The fact that $\zeta_i$ is dominated by
  $\omega_k$ also follows from the transversality of $W_i$ and
  $p(B_i)$.
\qed

As the referee pointed out, the construction above does yield a
cobordism in the case $k = 1$, which is topologically just a round
2--handle (a 2--handle and a 3--handle). It is not clear what can be
proved with such a cobordism, however, because it will not contain a
sphere of square 0.

\subsubsection*{Thanks}
The author would like to thank John Etnyre and Andras Stipsicz for
helpful comments on an initial version of this paper.

\bibliographystyle{gtart}
\bibliography{link}

\end{document}